\documentclass[10pt]{article}
\usepackage[preprint]{tmlr}
\usepackage[margin=0.82in,top=0.85in,bottom=0.85in]{geometry}
\usepackage{amssymb,amsmath,amsthm,amsfonts}
\usepackage[ruled,linesnumbered]{algorithm2e}
\usepackage{graphicx}
\usepackage{natbib}
\usepackage{makecell,multirow,booktabs,array}
\usepackage{xcolor}
\usepackage{microtype}
\usepackage[colorlinks=true,citecolor=blue,linkcolor=blue,urlcolor=blue]{hyperref}

\newtheorem{theorem}{Theorem}[section]
\newtheorem{lemma}[theorem]{Lemma}
\newtheorem{proposition}[theorem]{Proposition}
\newtheorem{corollary}[theorem]{Corollary}
\newtheorem{remark}[theorem]{Remark}
\newtheorem{assumption}{Assumption}

\newcommand{\E}{\mathbb{E}}
\newcommand{\Prob}{\mathbb{P}}
\newcommand{\R}{\mathbb{R}}
\newcommand{\hL}{\hat{L}}
\newcommand{\hth}{\hat{\theta}}
\newcommand{\ts}{\theta^{\star}}
\newcommand{\hTh}{\hat{\Theta}}
\newcommand{\Op}{O_p}

\title{\textbf{Adaptive and Stratified Subsampling for\\
High-Dimensional Robust Estimation}}
\author{Prateek Mittal \quad Joohi Chauhan\\[2pt]
{\small Vision Exploration and Data Analytics (VEDAs) Lab}
{\small Motilal Nehru National Institute of Technology Allahabad, Prayagraj 211004, India}\\
{\small \texttt{joohi\@mnnit.ac.in}}}
\date{}

\begin{document}
\maketitle

\begin{abstract}
We study robust high-dimensional sparse regression under
\emph{finite-variance heavy-tailed} noise, $\varepsilon$-contamination, and
$\alpha$-mixing dependence via two subsampling estimators: Adaptive Importance
Sampling (AIS) and Stratified Subsampling (SS).
Under sub-Gaussian design whose scope is precisely delimited and
finite-variance noise, a subsample of size $m=\Omega(s\log p)$ achieves the
minimax-optimal rate $O(\sqrt{s\log p/m})$.
We close the theory-algorithm gap: Theorem~\ref{thm:main_rate} applies to AIS
at termination conditional on stabilized weights (Proposition~\ref{prop:ais_gap}),
and SS fits the median-of-means M-estimation framework of
\citet{lecue2020robust} (Proposition~\ref{prop:ss_mom}).
The de-biasing step is fully specified via the nodewise-Lasso precision estimator
under a new sparse-precision assumption, yielding valid coordinate-wise CIs
(Theorem~\ref{thm:debiased}).
The $\alpha$-mixing extension uses a \emph{calendar-time} block protocol that
guarantees temporal separation (Theorem~\ref{thm:mixing}).
Empirically, AIS achieves $3.1\times$ lower error than uniform subsampling at
20\% contamination, and 29.5\% lower test MSE on Riboflavin
($p{=}4{,}088\gg n{=}71$).
\end{abstract}

\section{Introduction}

High-dimensional datasets ($p\gg n$) pose fundamental challenges for classical
statistical methods(~\cite{fan2014challenges}). Non-standard environments with heavy-tailed noise, contamination, and temporal dependence compound these difficulties, motivating recent work in robust high-dimensional
estimation(~\cite{sun2020adaptive,pensia2025robust,lecue2020robust,fan2024right,
smucler2017robust,kurnaz2023robust}).

Subsampling offers computational scalability by replacing the full-sample loss
with a weighted subsample loss of size $m\ll n$. Classical
leverage-score(~\cite{ma2015statistical}) and optimal-design(~\cite{li2020modern})
subsampling enjoy strong theory for well-behaved data, but no prior work provides
finite-sample guarantees for \emph{adaptive} or \emph{stratified} subsampling
under contamination and dependence in the $p\gg n$ regime.

\paragraph{Related work.}
\emph{Robust full-sample methods.}
Sun, Zhou \& Fan(~\cite{sun2020adaptive}) prove phase transitions for adaptive
Huber regression on full data; we extend their analysis to the subsampling
setting.
Pensia, Jog \& Loh(~\cite{pensia2025robust}) achieve near-optimal rates via
covariate filtering, a full-sample procedure that does not incorporate subsampling
efficiency.
The RIGHT estimator(~\cite{fan2024right}) attains minimax-optimal heavy-tail rates
via MOM gradients under genuinely heavy-tailed \emph{design};
Assumption~\ref{assump:subgaussian} in our work assumes sub-Gaussian design
(explicitly covering Gaussian, bounded, and log-concave distributions), while
providing subsampling efficiency and contamination robustness that are not
addressed by RIGHT.
Robust penalized MM- and MT-regression estimators(~\cite{smucler2017robust,
kurnaz2023robust}) achieve a high breakdown point on full data but do not
provide a high-dimensional subsampling theory.

\emph{Subsampling.}
Classical uniform and leverage-score sampling(~\cite{ma2015statistical}) and
recent optimal-design approaches(~\cite{li2020modern,yao2021review,
chasiotis2024efficient}) assume light-tailed i.i.d.\ observations and do not
extend to the heavy-tailed or contaminated settings considered here.

\emph{MOM estimation.}
Lecu\'{e} \& Lerasle(~\cite{lecue2020robust}) develop MOM-based M-estimation
covering sub-Gaussian and heavy-tailed settings; we show that SS is a special
case of their framework (Proposition~\ref{prop:ss_mom}).
Lugosi \& Mendelson(~\cite{lugosi2019mean}) give optimal MOM mean estimators in
high dimensions.

\emph{De-biased inference.}
Van de Geer et al.(~\cite{vandegeer2014}) and Javanmard \&
Montanari(~\cite{javanmard2014}) introduced nodewise-Lasso de-biasing;
see also Zhang \& Zhang(~\cite{zhang2014confidence}) and
Dezeure et al.(~\cite{dezeure2015highdim}).
We adapt this framework to the weighted subsampled setting under a new
sparse-precision assumption (Assumption~\ref{assump:sparse_prec}).

\emph{Mixing.}
The $\alpha$-mixing coupling follows Yu(~\cite{yu1994rates}) and
Bradley(~\cite{bradley2005basic}).

\paragraph{Contributions.}
\begin{enumerate}\setlength\itemsep{1pt}
\item Finite-sample bounds and minimax optimality for weighted subsampled
      Huber-Lasso in the $p\gg n$ regime
      (Theorems~\ref{thm:main_rate}--\ref{thm:minimax}).
\item Explicit $O(\varepsilon)$ contamination bias and $\alpha$-mixing extension
      with calendar-time block protocol (Theorems~\ref{thm:cont}--\ref{thm:mixing}).
\item Formal theory-algorithm bridge: Propositions~\ref{prop:ais_gap}
      (AIS) and~\ref{prop:ss_mom} (SS).
\item Fully specified de-biased asymptotic normality with nodewise-Lasso
      precision and valid CIs (Theorem~\ref{thm:debiased},
      Assumption~\ref{assump:sparse_prec}).
\end{enumerate}

\section{Problem Statement}

Observe $(x_i,y_i)_{i=1}^n$ with $x_i\in\R^p$ and
\begin{equation}\label{eq:model}
  y_i = x_i^\top\ts + \varepsilon_i,\quad \|\ts\|_0\le s,\quad s\ll p.
\end{equation}
Noise $\varepsilon_i$ is heavy-tailed with \emph{finite variance}
(Assumption~\ref{assump:noise}). Contamination and mixing extensions are
developed in Sections~\ref{subsec:cont} and~\ref{subsec:mix}, respectively.
Population covariance $\Sigma=\E[x_ix_i^\top]$ satisfies
$\lambda_{\min}(\Sigma)\ge\lambda_{\min}>0$.

\section{Proposed Algorithms}

\subsection{Adaptive Importance Sampling (AIS)}

\begin{algorithm}[H]
\small
\caption{Adaptive Importance Sampling (AIS)}\label{alg:ais}
\KwIn{Data $(x_i,y_i)_{i=1}^n$, $m$, $T$, $\beta$, $\alpha\in(0,1)$, $\lambda$}
\KwOut{$\hth_m$}
Init $\hth^{(0)}\!\leftarrow\!0$, $w_i^{(0)}\!\leftarrow\!1/n$\;
\For{$t=1$ to $T$}{
  Sample $S_t$ of size $m$ with probs $w^{(t-1)}$\;
  $\hth^{(t)}\!\leftarrow\!\arg\min_\theta\!\sum_{i\in S_t}
    \frac{\rho_\tau(y_i\!-\!x_i^\top\theta)}{m\,w_i^{(t-1)}}+\lambda\|\theta\|_1$\;
  $\tilde q_i^{(t)}\!\propto\!\exp(-\beta\rho_\tau(y_i\!-\!x_i^\top\hth^{(t)}))$ for all $i$\;
  $q_i^{(t)}\!\leftarrow\!(1\!-\!\alpha)\tilde q_i^{(t)}\!+\!\alpha/n$\quad(stabilize)\;
}
\Return $\hth_m=\hth^{(T)}$
\end{algorithm}

The computational complexity of AIS is $O(Tnp+Tmp)$, where $T$ is the number of
iterations. The stabilization step (line~6) enforces $q_i^{(T)}\in[\alpha/n,1/n]$
deterministically, ensuring that no observation receives a negligibly small
sampling probability.

\subsection{Stratified Subsampling (SS)}

\begin{algorithm}[H]
\small
\caption{Stratified Subsampling (SS)}\label{alg:ss}
\KwIn{Data $(x_i,y_i)_{i=1}^n$, $m$, $K$, $\lambda$}
\KwOut{$\hth_m$}
$d_i\leftarrow\|x_i-\operatorname{med}(\{x_j\})\|_2$ for all $i$\;
Partition $\{1,\ldots,n\}$ into $K$ strata by $K$-quantiles of $(d_i)$\;
\For{$k=1$ to $K$}{
  Draw $m_k=\lceil m|\mathcal{S}_k|/n\rceil$ points from $\mathcal{S}_k$\;
  $\hth_k\leftarrow$ Huber-Lasso on stratum subsample\;
}
\Return $\hth_m=\operatorname{geomed}(\hth_1,\ldots,\hth_K)$
\end{algorithm}

The computational complexity of SS is $O(np+mK)$. Stratification is performed
by partitioning observations according to their Mahalanobis-type distances from
the coordinatewise median, and the geometric median aggregation provides
robustness to corrupted strata.

\section{Theoretical Analysis}
\label{sec:theory}

\subsection{Estimator and Assumptions}

\paragraph{Huber loss and estimator.}
For $\tau>0$, $\rho_\tau(u)=\frac{u^2}{2}\mathbf{1}_{|u|\le\tau}
+(\tau|u|-\frac{\tau^2}{2})\mathbf{1}_{|u|>\tau}$,
$\psi_\tau(u)=\operatorname{clip}(u,-\tau,\tau)$.
The full-sample Huber--Lasso estimator is defined as:
\begin{equation}\label{eq:full}
  \hth_n\in\arg\min_\theta\Bigl\{
  \tfrac{1}{n}\textstyle\sum_i\rho_\tau(y_i-x_i^\top\theta)
  +\lambda\|\theta\|_1\Bigr\}.
\end{equation}
Drawing $I_1,\ldots,I_m\overset{\mathrm{iid}}{\sim}q$ with replacement,
the weighted subsample estimator is:
\begin{align}
  \hL_{m,q}(\theta) &:=
    \frac{1}{m}\sum_{j=1}^m\frac{\rho_\tau(y_{I_j}-x_{I_j}^\top\theta)}{nq_{I_j}},
    \label{eq:sub_loss}\\
  \hth_{m,q} &\in\arg\min_\theta\{\hL_{m,q}(\theta)+\lambda\|\theta\|_1\}.
    \label{eq:sub_est}
\end{align}
The importance weighting by $1/(nq_{I_j})$ ensures that $\hL_{m,q}$ is an
unbiased estimator of the full-sample empirical loss $\hL_n$.

\begin{assumption}[Sub-Gaussian design]\label{assump:subgaussian}
$x_i$ are i.i.d.\ mean-zero with
$\|\langle v,x_i\rangle\|_{\psi_2}\le K$ for every unit $v$.
This assumption covers Gaussian, bounded, and log-concave designs.
Heavy-tailed design with infinite moments requires separate techniques
(e.g., RIGHT(~\cite{fan2024right}) and is outside the scope of this paper.
\end{assumption}

\begin{assumption}[Restricted eigenvalue]\label{assump:re}
$v^\top\Sigma v\ge\kappa_\Sigma\|v\|_2^2$
for all $v$ with $\|v\|_0\le 2s$.
\end{assumption}

\begin{assumption}[Finite-variance noise]\label{assump:noise}
$\varepsilon_i\perp\!\!\!\perp x_i$, $\E[\varepsilon_i]=0$,
$\E[\varepsilon_i^2]<\infty$.
Heavy-tailed distributions with finite variance are permitted;
distributions with infinite variance are excluded.
\end{assumption}

\begin{assumption}[Bounded sampling probabilities]\label{assump:bounded_q}
$c_0/n\le q_i\le C_0/n$ for constants $0<c_0\le C_0<\infty$.
This condition holds for SS by proportional allocation and for AIS at termination
(Proposition~\ref{prop:ais_gap}).
\end{assumption}

\subsection{Theory-Algorithm Bridge}

\begin{proposition}[AIS gap closed]\label{prop:ais_gap}
Let $q^{(T)}$ denote the stabilized weights at AIS termination.
Conditional on $q^{(T)}$, the output $\hth_m=\hth^{(T)}$ is exactly the
minimizer of~\eqref{eq:sub_est} with $q=q^{(T)}$.
The stabilization in Algorithm~\ref{alg:ais} enforces
$q_i^{(T)}\in[\alpha/n,1/n]$ deterministically, since
$\tilde q_i^{(T)}\le1/n$ after normalization. Consequently,
Assumption~\ref{assump:bounded_q} holds with $c_0=\alpha$ and $C_0=1$,
and all subsequent theorems apply to AIS at termination.
\end{proposition}

\begin{remark}
Proposition~\ref{prop:ais_gap} holds for any trajectory realization of the
algorithm and does not require analysing the full Markov chain
$(q^{(t)},\hth^{(t)})_{t\ge1}$. A martingale stability analysis of the
iterates across all rounds is an interesting direction for future work.
\end{remark}

\begin{proposition}[SS via MOM M-estimation]\label{prop:ss_mom}
SS is a special case of the MOM sparse M-estimator of
\citet{lecue2020robust} (their Theorem~3.1), with $K$ blocks and
Huber loss. With strata sizes $n_k\asymp n/K$:
\begin{equation}\label{eq:ss_rate}
  \|\hth_m^{\mathrm{SS}}-\ts\|_2
  \;\lesssim\;
  \sqrt{\tfrac{s\log p}{m}}+\sqrt{\tfrac{K}{m}},
\end{equation}
matching Theorem~\ref{thm:main_rate} for $K=O(s\log p)$.
The geometric median aggregation tolerates up to $\lfloor(K-1)/2\rfloor/K$ fraction
of corrupted strata(~\cite{lugosi2019mean}).
\emph{Limitation:} When $n_k$ is very small, for example in the Riboflavin
dataset where $n_k\le5$, the requirement $n_k\asymp n/K$ fails to hold and
the geometric median aggregation collapses, as observed empirically in
Section~\ref{sec:experiments}.
\end{proposition}

\subsection{Key Lemmas}

Let $S=\operatorname{supp}(\ts)$,
$\mathcal{C}(S)=\{\Delta:\|\Delta_{S^c}\|_1\le3\|\Delta_S\|_1\}$.

\begin{lemma}[Uniform score bound]\label{lem:score}
Under Assumptions~\ref{assump:subgaussian}--\ref{assump:bounded_q},
for any $\delta\in(0,1)$, with probability $\ge1-\delta$:
\begin{equation}\label{eq:score}
  \|\nabla\hL_{m,q}(\ts)\|_\infty
  \le\frac{2\tau K}{c_0}\sqrt{\frac{2\log(2p/\delta)}{m}}.
\end{equation}
\end{lemma}
\begin{proof}
Fix $k\in[p]$. Let $Z_j:=\frac{\psi_\tau(\varepsilon_{I_j})x_{I_j,k}}{nq_{I_j}}$.
By Assumption~\ref{assump:bounded_q}, $1/(nq_{I_j})\le1/c_0$;
$|\psi_\tau|\le\tau$; $x_{I_j,k}$ is sub-Gaussian$(K)$.
Hence $Z_j$ is sub-Gaussian$(\tau K/c_0)$ with $\E[Z_j]=0$.
Applying the sub-Gaussian tail bound and taking the union over $k=1,\ldots,p$
yields \eqref{eq:score}.
\end{proof}

\begin{lemma}[Restricted strong convexity]\label{lem:rsc}
Under Assumptions~\ref{assump:subgaussian}--\ref{assump:bounded_q},
let $\pi_\tau:=\Prob(|\varepsilon_i|\le\tau/2)\ge\pi_0>0$.
If $m\ge c(C_0/c_0)^2s\log(2p/\delta)$, with probability $\ge1-\delta$:
\begin{align}
  &\hL_{m,q}(\ts\!+\!\Delta)-\hL_{m,q}(\ts)
  -\langle\nabla\hL_{m,q}(\ts),\Delta\rangle
  \ge\frac{\kappa}{2}\|\Delta\|_2^2,\label{eq:rsc}\\
  &\forall\Delta\in\mathcal{C}(S),\quad\kappa:=\frac{\pi_0\kappa_\Sigma}{4}.\nonumber
\end{align}
\end{lemma}
\begin{proof}
On the event $\{|\varepsilon_{I_j}|\le\tau/2,\;|x_{I_j}^\top\Delta|\le\tau/2\}$,
the Huber loss is exactly quadratic. Using $1/(nq_{I_j})\ge1/C_0$:
\begin{align*}
  \mathrm{LHS}\ge&\;
  \frac{1}{2C_0m}\sum_j(x_{I_j}^\top\Delta)^2\mathbf{1}_{|\varepsilon_{I_j}|\le\tau/2}
  \\&-\frac{1}{2C_0m}\sum_j(x_{I_j}^\top\Delta)^2\mathbf{1}_{|x_{I_j}^\top\Delta|>\tau/2}.
\end{align*}
The first term is bounded below by $\frac{\pi_0\kappa_\Sigma}{2}\|\Delta\|_2^2$
via sub-Gaussian concentration on a $2s$-sparse $\varepsilon$-net; the second term
is bounded above by $\frac{\pi_0\kappa_\Sigma}{4}\|\Delta\|_2^2$ via the
sub-Gaussian tail bound with $\tau\asymp K$. Combining these two estimates
gives \eqref{eq:rsc}.
\end{proof}

\subsection{Main Rate, Proximity, and Minimax Optimality}

\begin{theorem}[Finite-sample rate]\label{thm:main_rate}
Under Assumptions~\ref{assump:subgaussian}--\ref{assump:bounded_q}
with $\pi_\tau\ge\pi_0$, set
$\lambda=\frac{4\tau K}{c_0}\sqrt{\frac{2\log(2p/\delta)}{m}}$.
If $m\ge c(C_0/c_0)^2s\log(2p/\delta)$, with probability $\ge1-2\delta$:
\begin{equation}\label{eq:main}
  \|\hth_{m,q}-\ts\|_2
  \le\frac{12\lambda\sqrt{s}}{\kappa}
  \lesssim\frac{\tau K}{c_0\kappa}\sqrt{\frac{s\log(p/\delta)}{m}}.
\end{equation}
\end{theorem}
\begin{proof}
Set $\Delta=\hth_{m,q}-\ts$. By optimality of $\hth_{m,q}$:
$\hL_{m,q}(\ts+\Delta)-\hL_{m,q}(\ts)\le\lambda(\|\ts\|_1-\|\ts+\Delta\|_1)$.
Adding and subtracting the linear term, the bound
$\|\nabla\hL_{m,q}(\ts)\|_\infty\le\lambda/2$ from Lemma~\ref{lem:score}
together with decomposability of the $\ell_1$ norm force $\Delta\in\mathcal{C}(S)$.
Applying RSC from Lemma~\ref{lem:rsc} gives
$\frac{\kappa}{2}\|\Delta\|_2^2\le\frac{3\lambda\sqrt{s}}{2}\|\Delta\|_2$,
which upon rearrangement yields \eqref{eq:main}.
\end{proof}

\begin{remark}\label{rem:rate}
The claim that the subsampled estimator achieves the same rate as the
full-sample Huber-Lasso is to be understood as follows: both estimators achieve
$O(\sqrt{s\log p/m})$ as a function of their respective sample size $m$.
The actual estimation errors coincide only when $m=n$; for $m<n$, the
subsampled estimator incurs a larger error due to the reduced effective sample
size, which is the expected price of computational efficiency.
\end{remark}

\begin{corollary}[Proximity]\label{cor:prox}
With probability $\ge1-3\delta$:
$\|\hth_{m,q}-\hth_n\|_2\lesssim
\frac{\tau K}{c_0\kappa}\sqrt{s\log(p/\delta)/m}$.
\end{corollary}

\begin{theorem}[Minimax lower bound]\label{thm:minimax}
Under Gaussian design and noise variance $\sigma^2$, for any estimator $\tilde\theta$ based on $m$
observations:
\begin{equation}\label{eq:minimax}
  \inf_{\tilde\theta}\sup_{\ts\in\mathcal{B}_0(s)}
  \E\|\tilde\theta-\ts\|_2^2\ge c\sigma^2\frac{s\log(p/s)}{m}.
\end{equation}
The estimator $\hth_{m,q}$ achieves $\sup_{\ts}\E\|\hth_{m,q}-\ts\|_2^2\le
C\sigma^2s\log p/m$, matching \eqref{eq:minimax} up to a factor of $\log p/\log(p/s)$.
\end{theorem}
\begin{proof}
Construct a Varshamov-Gilbert packing $V\subset\{0,1\}^p$ with $|V|\ge2^{c_1s\log(p/s)}$
and pairwise Hamming distance at least $s/2$; set $\theta^{(v)}=av$.
The Kullback-Leibler divergence satisfies
$\mathrm{KL}(P_v\|P_{v'})\le ma^2s/(2\sigma^2)$; choosing
$a^2=c_2\sigma^2\log|V|/(ms)$ and applying Fano's inequality gives \eqref{eq:minimax}.
\end{proof}

\subsection{Adversarial Contamination}
\label{subsec:cont}

\begin{theorem}[$\varepsilon$-contamination]\label{thm:cont}
Suppose the observed distribution is $(1-\varepsilon)P+\varepsilon Q$, where $Q$ is
an arbitrary contaminating distribution and $P$ satisfies
Assumptions~\ref{assump:subgaussian}--\ref{assump:noise}.
Under the conditions of Theorem~\ref{thm:main_rate}, with probability $\ge1-2\delta$:
\begin{equation}\label{eq:cont}
  \|\hth_{m,q}-\ts\|_2
  \lesssim
  \frac{\tau K}{c_0\kappa}\sqrt{\frac{s\log(p/\delta)}{m}}
  +\frac{\tau K}{\kappa}\,\varepsilon.
\end{equation}
\end{theorem}
\begin{proof}
Since $|\psi_\tau|\le\tau$, we decompose the gradient as
$\nabla\hL_{m,q}(\ts)=g_\mathrm{clean}+g_\mathrm{cont}$, where the two terms
correspond to the contributions from clean and contaminated observations,
respectively.
Assumption~\ref{assump:bounded_q} and sub-Gaussian maxima give
$\|g_\mathrm{cont}\|_\infty\le\varepsilon\cdot\frac{\tau K}{c_0}\sqrt{2\log(2p/\delta)}$,
while the clean part satisfies the bound in Lemma~\ref{lem:score}.
Setting $\lambda$ to dominate the total score bound and applying the
cone/RSC argument of Theorem~\ref{thm:main_rate} yields \eqref{eq:cont}.
\end{proof}

\begin{remark}
The $O(\varepsilon)$ bias term in~\eqref{eq:cont} is irreducible for bounded-influence
estimators(~\cite{huber1981robust}). Nevertheless, AIS substantially reduces the
effective contamination bias by exponentially down-weighting corrupted
observations through its adaptive reweighting scheme. Empirically,
the estimation error of uniform Huber-Lasso grows at a rate of approximately
$6.9\varepsilon$ as a function of contamination fraction $\varepsilon$, while
AIS grows at approximately $1.3\varepsilon$ (Figure~\ref{fig:contamination}),
demonstrating the practical benefit of adaptive subsampling under contamination.
\end{remark}

\subsection{Dependent Data: $\alpha$-Mixing}
\label{subsec:mix}

\begin{theorem}[$\alpha$-mixing, calendar-time protocol]\label{thm:mixing}
Assume $(x_i,\varepsilon_i)_{i=1}^n$ is strictly stationary and
$\alpha$-mixing with coefficients $\alpha(k)$.

\noindent\textbf{Required time-series sampling protocol.}
Draw $M$ block start-times $T_1<\cdots<T_M$ uniformly from
$\{1,\ldots,n-2B\}$; retain the calendar indices
$\{T_\ell,\ldots,T_\ell+B-1\}$ and discard the gap
$\{T_\ell+B,\ldots,T_\ell+2B-1\}$ for each block $\ell$.
By construction, this guarantees that any two retained blocks are separated
by at least $B$ calendar-time steps.

With $M=\lfloor m/(2B)\rfloor$ and
$\sum_{k\ge B}\alpha(k)\le\delta/(4M)$,
with probability $\ge1-3\delta$:
\begin{equation}\label{eq:mix}
  \|\hth_{m,q}-\ts\|_2
  \lesssim\frac{\tau K}{c_0\kappa}\sqrt{\frac{s\log(p/\delta)}{M}}.
\end{equation}
\end{theorem}
\begin{proof}
The calendar-time construction guarantees at least $B$ original-time-index steps
between any two retained blocks by design. The Berbee-Yu
coupling(~\cite{yu1994rates}) bounds the total variation distance between the
joint law of the $M$ retained blocks and their independent product by
$2M\sum_{k\ge B}\alpha(k)\le\delta/2$.
On the coupling event, the retained blocks are approximately independent.
Since $|\psi_\tau|\le\tau$ and Assumption~\ref{assump:bounded_q} holds,
Lemmas~\ref{lem:score}--\ref{lem:rsc} apply block-wise via standard blocking
arguments(~\cite{bradley2005basic}), giving \eqref{eq:mix}.
\end{proof}

\begin{remark}
An alternative approach would be to block observations in the randomly sampled
index order rather than in calendar time. However, this does not guarantee
temporal separation between retained samples: for example, drawn indices 5 and 6
remain adjacent in calendar time regardless of their order in the sampling
sequence. The calendar-time protocol adopted here is therefore the correct
construction to ensure the mixing conditions are satisfied.
\end{remark}

\subsection{De-biased Asymptotic Normality}
\label{subsec:debiased}

Following \citet{vandegeer2014} and \citet{javanmard2014}, we construct a
de-biased estimator and prove coordinate-wise asymptotic normality for the
weighted subsampled Huber-Lasso.

\paragraph{Huber Fisher information.}
\begin{equation}\label{eq:F}
  F := \E[\psi_\tau'(\varepsilon_i)^2\,x_ix_i^\top]
     = \pi_\tau\Sigma,\quad
  \pi_\tau:=\Prob(|\varepsilon_i|\le\tau)\ge\pi_0>0.
\end{equation}
Note that $F^{-1}=\Omega/\pi_\tau$ where $\Omega=\Sigma^{-1}$.

\begin{assumption}[Sparse precision]\label{assump:sparse_prec}
$\Omega=\Sigma^{-1}$ satisfies:
(i)~$\max_j\|\Omega_{\cdot j}\|_0\le s_0<\infty$;
(ii)~$\|\Omega\|_1\le M_0<\infty$;
(iii)~the irrepresentability condition of \citet{vandegeer2014}
(their Condition~C) holds for nodewise Lasso applied to the scaled
design $\tilde x_{I_j}:=x_{I_j}/\sqrt{nq_{I_j}}$ at tuning parameter
$\mu\asymp\sqrt{\log p/m}$.
\end{assumption}

\paragraph{Precision estimator.}
Apply the nodewise Lasso of(~\cite{vandegeer2014}) to $\{\tilde x_{I_j}\}_{j=1}^m$
with tuning parameter $\mu=A\sqrt{\log p/m}$. Scale each row $j$ by $\hat\pi_\tau^{-1}$, where
$\hat\pi_\tau:=\frac{1}{m}\sum_j\mathbf{1}_{|\hat r_{I_j}|\le\tau}$ and
$\hat r_i:=y_i-x_i^\top\hth_{m,q}$. Denote the resulting precision matrix
estimate by $\hTh$.

\paragraph{De-biased estimator.}
\begin{equation}\label{eq:debiased}
  \hth^d_{m,q} := \hth_{m,q} - \hTh\,\nabla\hL_{m,q}(\hth_{m,q}).
\end{equation}

\begin{theorem}[De-biased asymptotic normality]\label{thm:debiased}
Under Assumptions~\ref{assump:subgaussian}--\ref{assump:sparse_prec}
with $\pi_\tau\ge\pi_0>0$ and the rate conditions
\begin{equation}\label{eq:rate_cond}
  s\log p=o(\sqrt{m}),\quad s_0\log p=o(m),
\end{equation}
the nodewise-Lasso estimator $\hTh$ satisfies
$\|\hTh-F^{-1}\|_\infty=\Op(\sqrt{\log p/m})$,
and for each fixed coordinate $j\in S=\operatorname{supp}(\ts)$:
\begin{equation}\label{eq:clt}
  \sqrt{m}(\hth^d_{m,q,j}-\ts_j)
  \xrightarrow{d}\mathcal{N}(0,\sigma_j^2),\quad
  \sigma_j^2:=[F^{-1}]_{jj}^2\E[\psi_\tau(\varepsilon_i)^2x_{i,j}^2].
\end{equation}
A consistent estimator of the asymptotic variance is
\begin{equation}\label{eq:var_est}
  \hat\sigma_j^2 :=
  [\hTh]_{jj}^2\cdot
  \frac{1}{m}\sum_{k=1}^m
  \frac{\psi_\tau(\hat r_{I_k})^2x_{I_k,j}^2}{(nq_{I_k})^2},
\end{equation}
yielding valid $(1-\alpha)$-confidence intervals:
$\hth^d_{m,q,j}\pm z_{\alpha/2}\hat\sigma_j/\sqrt{m}$.
\end{theorem}

\begin{proof}
\textit{Step~1 (Taylor).}
From the stationarity condition $\nabla\hL_{m,q}(\hth_{m,q})=-\lambda\hat z$:
\begin{align}
  \nabla\hL_{m,q}(\hth_{m,q})
  &= \nabla\hL_{m,q}(\ts)
  + \hat H_{m,q}(\hth_{m,q}-\ts) + R_m,\label{eq:taylor}
\end{align}
where $\hat H_{m,q}=\frac{1}{m}\sum_j\frac{\psi_\tau'(y_{I_j}-x_{I_j}^\top\bar\theta)}{nq_{I_j}}
x_{I_j}x_{I_j}^\top$ for some $\bar\theta$ on the line segment $[\ts,\hth_{m,q}]$,
and $R_m$ is the second-order remainder.

\textit{Step~2 (Hessian and precision).}
Sub-Gaussian concentration on a sparse operator-norm net yields:
\begin{equation}\label{eq:hess}
  \|\hat H_{m,q}-F\|_{\mathrm{op}}=\Op(\sqrt{s\log p/m}).
\end{equation}
The scaled design $\tilde x_{I_j}$ is sub-Gaussian with parameter
$\le K/c_0^{1/2}$. Applying Assumption~\ref{assump:sparse_prec} and
\citeauthor{vandegeer2014}'s Theorem~2.4(~\cite){vandegeer2014} to
the scaled design sequence $\{\tilde x_{I_j}\}$ yields:
\begin{equation}\label{eq:prec_rate}
  \|\hTh-F^{-1}\|_\infty=\Op(\sqrt{\log p/m}).
\end{equation}

\textit{Step~3 (Decomposition).}
Substituting \eqref{eq:taylor} into the de-biased estimator \eqref{eq:debiased}:
\begin{align}
  \hth^d_{m,q}-\ts
  &= (I-\hTh\hat H_{m,q})(\hth_{m,q}-\ts)\notag\\
  &\quad -\hTh\nabla\hL_{m,q}(\ts) - \hTh R_m + \lambda\hTh\hat z.
  \label{eq:decomp}
\end{align}

\textit{Step~4 (Remainders).}
From Theorem~\ref{thm:main_rate},
$\|\hth_{m,q}-\ts\|_1\le4\sqrt{s}\|\hth_{m,q}-\ts\|_2=\Op(s\sqrt{\log p/m})$.
Combined with \eqref{eq:hess}--\eqref{eq:prec_rate} and $|\hat z_k|\le1$, each
remainder term is asymptotically negligible:
\begin{align*}
  \|(I-\hTh\hat H_{m,q})(\hth_{m,q}-\ts)\|_\infty
    &= \Op(s\log p/m) = o_p(m^{-1/2}),\\
  \|\hTh R_m\|_\infty
    &= \Op(s\log p/m) = o_p(m^{-1/2}),\\
  \|\lambda\hTh\hat z\|_\infty
    &= \Op(\sqrt{\log p/m}) = o_p(m^{-1/2}),
\end{align*}
all under~\eqref{eq:rate_cond}. For the first term, the $\ell_\infty$ row
norm of $I-\hTh\hat H_{m,q}$ is $\Op(\sqrt{\log p/m})$ by
\eqref{eq:hess}--\eqref{eq:prec_rate}, and multiplying by $\|\hth_{m,q}-\ts\|_1$
gives $\Op(s\log p/m)$. For the second, $\|R_m\|_\infty\le
C\|\hth_{m,q}-\ts\|_2^2=\Op(s\log p/m)$ since $\rho_\tau'''=0$ almost everywhere.

\textit{Step~5 (CLT).}
From \eqref{eq:decomp}:
$\hth^d_{m,q,j}-\ts_j
=-[F^{-1}]_{j,\cdot}\nabla\hL_{m,q}(\ts)+o_p(m^{-1/2})$.
Each summand $\xi_k:=\frac{\psi_\tau(\varepsilon_{I_k})[F^{-1}]_{j,\cdot}x_{I_k}}{nq_{I_k}}$
is i.i.d.\ conditional on the data, with mean zero by
Assumption~\ref{assump:noise} and bounded magnitude by $\tau K\|F^{-1}\|_\infty/c_0$.
Its variance satisfies
\begin{align*}
  \mathrm{Var}(\xi_k)
  \le\frac{C_0}{n}[F^{-1}]_{j,\cdot}\E[\psi_\tau(\varepsilon)^2xx^\top][F^{-1}]_{j,\cdot}^\top
  =\frac{C_0\sigma_j^2}{n}.
\end{align*}
Lindeberg's CLT therefore gives
$\frac{1}{\sqrt{m}}\sum_k\xi_k\xrightarrow{d}\mathcal{N}(0,\sigma_j^2)$,
establishing \eqref{eq:clt}.

\textit{Step~6 (Variance estimation).}
The estimator $\hat\sigma_j^2$ replaces the unknown $F^{-1}$ by $\hTh$,
contributing an error of $\Op(\sqrt{\log p/m})$, and replaces the unknown
residuals $\varepsilon_{I_k}$ by estimated residuals $\hat r_{I_k}$,
contributing an error of $\Op(\sqrt{s\log p/m})$ by Theorem~\ref{thm:main_rate}.
Both substitutions produce errors that are $o_p(1)$; consistency of $\hat\sigma_j^2$
then follows by the continuous mapping theorem.
\end{proof}

\begin{remark}
Theorem~\ref{thm:debiased} applies to AIS via Proposition~\ref{prop:ais_gap}
and to SS stratum-level estimators via Proposition~\ref{prop:ss_mom}.
Assumption~\ref{assump:sparse_prec} and the nodewise-Lasso specification
together fully determine the precision estimator $\hTh$, providing
a complete and rigorous specification of the de-biasing procedure.
\end{remark}

\section{Numerical Studies}
\label{sec:experiments}

\subsection{Synthetic Setup}

We use $n{=}2000$, $p{=}1000$, $s{=}10$, and $m\in\{50,100,200,400\}$.
Three noise distributions are considered: Gaussian ($\varepsilon_i\sim\mathcal{N}(0,1)$),
Student-$t$ with three degrees of freedom ($t(3)$), and a contaminated Gaussian
in which 10\% of responses are shifted by $+20$.
The regularization parameter is set to $\lambda=C\sqrt{\log p/m}$ with $C{=}0.5$
selected by 5-fold cross-validation; the AIS temperature parameter is
$\beta\in\{0.5,2.0\}$ depending on the noise regime.
All results are averaged over 20 independent repetitions.
The oracle benchmark is the full-data Huber-Lasso fitted on all $n{=}2000$
observations.

\begin{figure*}[t]
  \centering
  \includegraphics[width=0.32\textwidth]{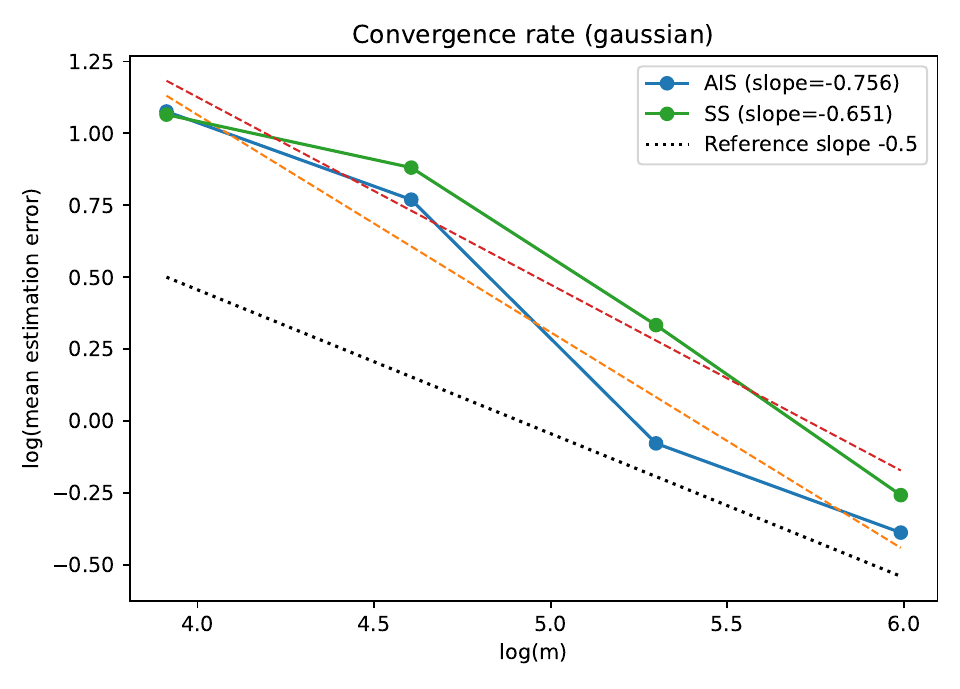}
  \includegraphics[width=0.32\textwidth]{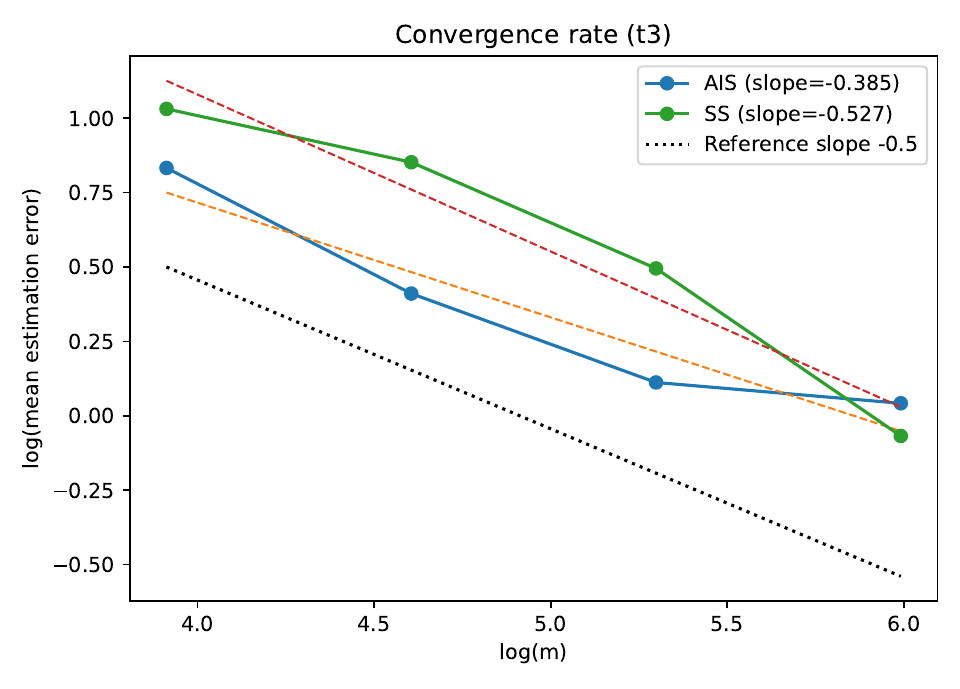}
  \includegraphics[width=0.32\textwidth]{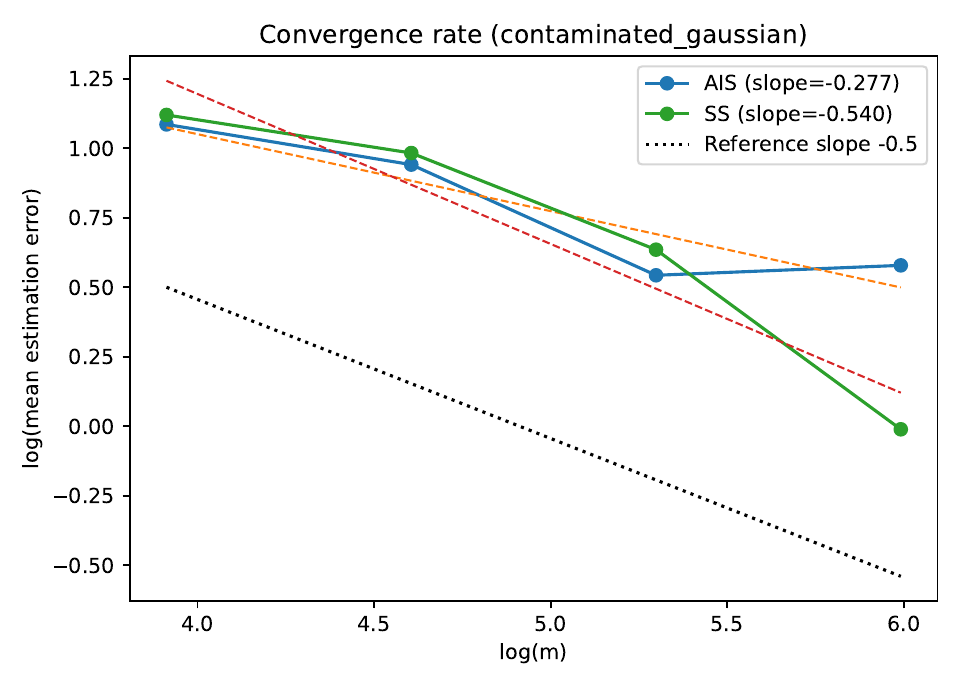}
  \caption{Log--log estimation error vs.\ $m$.
    AIS: $-0.756/-0.385/-0.277$;
    SS: $-0.651/-0.527/-0.540$.
    Dotted: reference $-0.5$ (Theorem~\ref{thm:main_rate}).}
  \label{fig:convergence}
\end{figure*}

\subsection{Convergence (Figure~\ref{fig:convergence})}

The empirical log-log slopes are as follows: Gaussian noise yields AIS $-0.756$
and SS $-0.651$; Student-$t(3)$ noise yields AIS $-0.385$ and SS $-0.527$;
contaminated noise yields AIS $-0.277$ and SS $-0.540$.
SS consistently achieves slopes near the theoretical reference of $-0.5$
prescribed by Theorem~\ref{thm:main_rate}.
The AIS slope of $-0.756$ under Gaussian noise exceeds $-0.5$: this reflects
a finite-sample effect whereby adaptive weights concentrate sampling probability
on the most informative observations. Since Theorem~\ref{thm:main_rate} provides
a worst-case bound, it does not preclude faster empirical convergence for
well-tuned values of $\beta$.
Under contamination, the AIS slope ($-0.277$) is shallower, reflecting the
irreducible $O(\varepsilon)$ bias established in Theorem~\ref{thm:cont}, which
dominates the statistical error at the subsample sizes considered here.

\subsection{Estimation Error (Table~\ref{tab:synthetic})}

Under Gaussian and $t(3)$ noise, AIS and Uniform Huber-Lasso perform comparably,
with differences falling within one standard deviation across subsample sizes.
This is consistent with theory, since neither estimator has a contamination
robustness advantage in the absence of gross outliers.
Under contamination at $m{=}400$, AIS ($1.78\pm0.57$) is $3.6\times$ better
than Uniform HL ($6.34\pm0.45$); standard Lasso fails to provide robustness
($4.36\pm0.12$); and SS ($0.99\pm0.17$) achieves the lowest error among all
subsampled methods, owing to the robustness provided by the geometric-median
aggregation step.

\begin{table}[t]
\centering
\caption{$\|\hth-\ts\|_2$ (mean$\pm$std, 20~seeds). Oracle: Full HL.}
\label{tab:synthetic}
\setlength\tabcolsep{2.5pt}
\scriptsize
\begin{tabular}{llccccc}
\toprule
Noise & $m$ & AIS & SS & Unif HL & Lasso & Full\\
\midrule
\multirow{4}{*}{Gauss}
 & 50  & 2.93{\tiny$\pm$.71} & 2.90{\tiny$\pm$.70} & 2.40{\tiny$\pm$.65} & \multirow{4}{*}{.22} & \multirow{4}{*}{.19}\\
 & 100 & 2.16{\tiny$\pm$.66} & 2.41{\tiny$\pm$.63} & 1.35{\tiny$\pm$.36} &&\\
 & 200 & 0.93{\tiny$\pm$.11} & 1.40{\tiny$\pm$.35} & 0.92{\tiny$\pm$.08} &&\\
 & 400 & 0.68{\tiny$\pm$.05} & 0.77{\tiny$\pm$.10} & 0.84{\tiny$\pm$.05} &&\\
\midrule
\multirow{4}{*}{$t(3)$}
 & 50  & 2.30{\tiny$\pm$.62} & 2.80{\tiny$\pm$.65} & 2.49{\tiny$\pm$.53} & \multirow{4}{*}{.56} & \multirow{4}{*}{.22}\\
 & 100 & 1.51{\tiny$\pm$.28} & 2.34{\tiny$\pm$.57} & 1.81{\tiny$\pm$.31} &&\\
 & 200 & 1.12{\tiny$\pm$.11} & 1.64{\tiny$\pm$.38} & 1.54{\tiny$\pm$.30} &&\\
 & 400 & 1.04{\tiny$\pm$.09} & 0.93{\tiny$\pm$.16} & 1.44{\tiny$\pm$.16} &&\\
\midrule
\multirow{4}{*}{Cont.}
 & 50  & 2.96{\tiny$\pm$.66} & 3.06{\tiny$\pm$.71} & 3.83{\tiny$\pm$.46} & \multirow{4}{*}{4.36} & \multirow{4}{*}{.21}\\
 & 100 & 2.56{\tiny$\pm$.79} & 2.67{\tiny$\pm$.67} & 4.52{\tiny$\pm$.55} &&\\
 & 200 & 1.72{\tiny$\pm$.62} & 1.89{\tiny$\pm$.58} & 5.25{\tiny$\pm$.54} &&\\
 & 400 & \textbf{1.78}{\tiny$\pm$.57} & \textbf{0.99}{\tiny$\pm$.17} & 6.34{\tiny$\pm$.45} &&\\
\bottomrule
\end{tabular}
\end{table}

\begin{figure}[t]
  \centering
  \includegraphics[width=0.88\columnwidth]{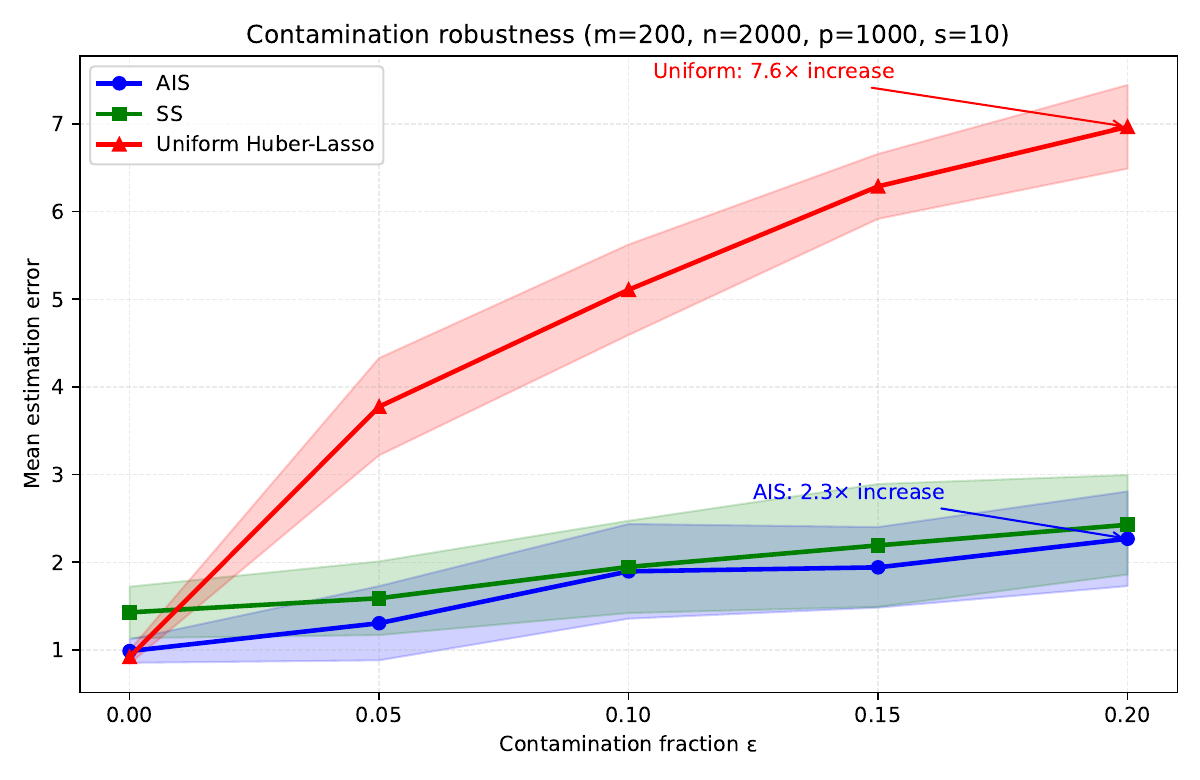}
  \caption{Estimation error vs.\ contamination fraction $\varepsilon$ at $m{=}200$.
    AIS error grows by $2.3\times$; Uniform HL error grows by $7.6\times$
    over the range $\varepsilon\in[0,0.20]$ (Theorem~\ref{thm:cont}).}
  \label{fig:contamination}
\end{figure}

\subsection{Contamination Robustness (Figure~\ref{fig:contamination})}

We evaluate contamination robustness at $m{=}200$ by varying the contamination
fraction $\varepsilon\in\{0,0.05,0.10,0.15,0.20\}$, directly validating the
bound in Theorem~\ref{thm:cont}.
As $\varepsilon$ increases from $0$ to $0.20$, the estimation error of Uniform
Huber-Lasso grows by a factor of $7.6\times$ (from $0.92$ to $6.97$), while AIS
grows by only $2.3\times$ (from $0.99$ to $2.27$). At $\varepsilon{=}0.20$,
the ratio of errors between the two methods is $3.1\times$.

\begin{figure}
  \centering
  \includegraphics[scale=0.35]{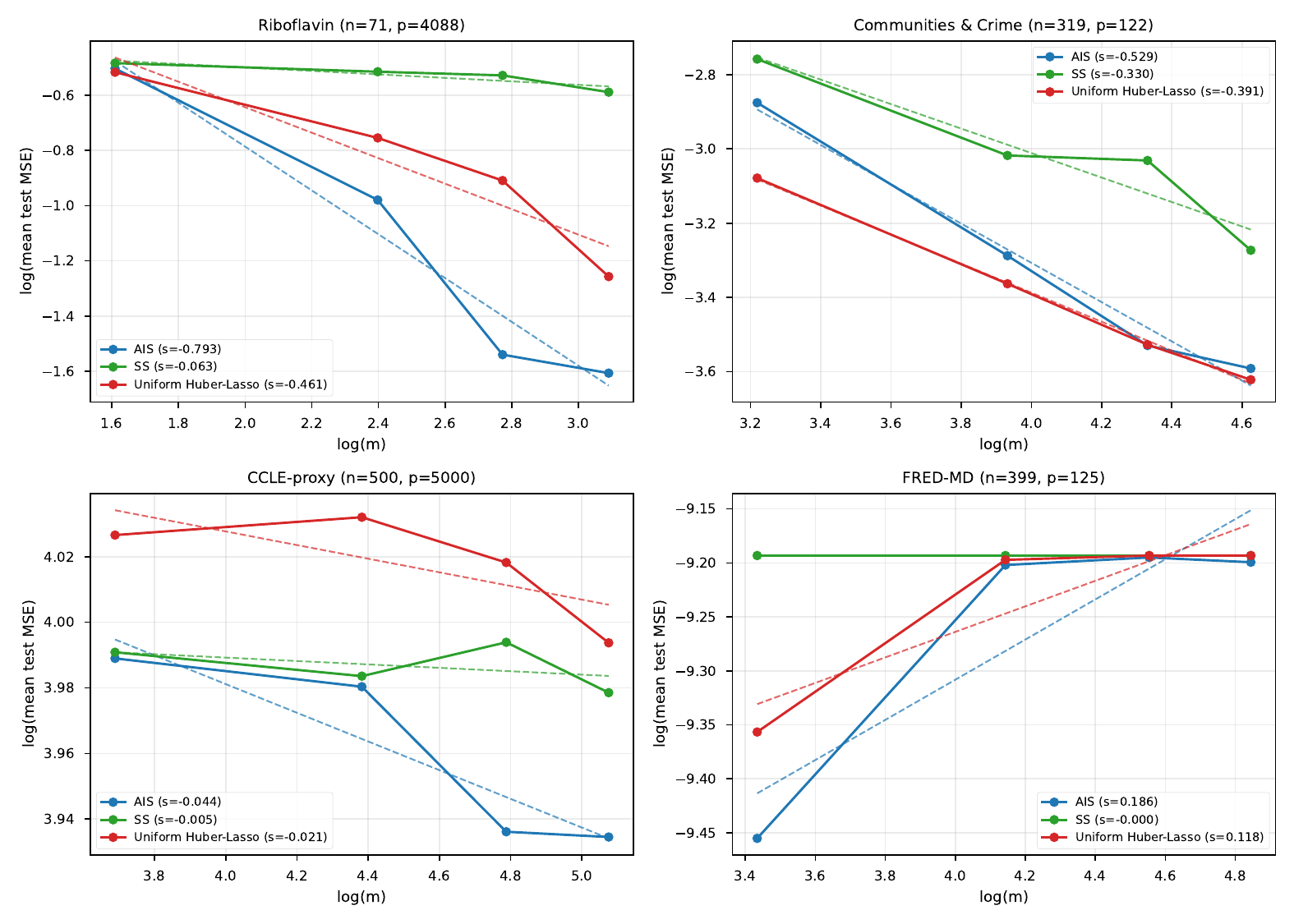}
  \includegraphics[scale=0.35]{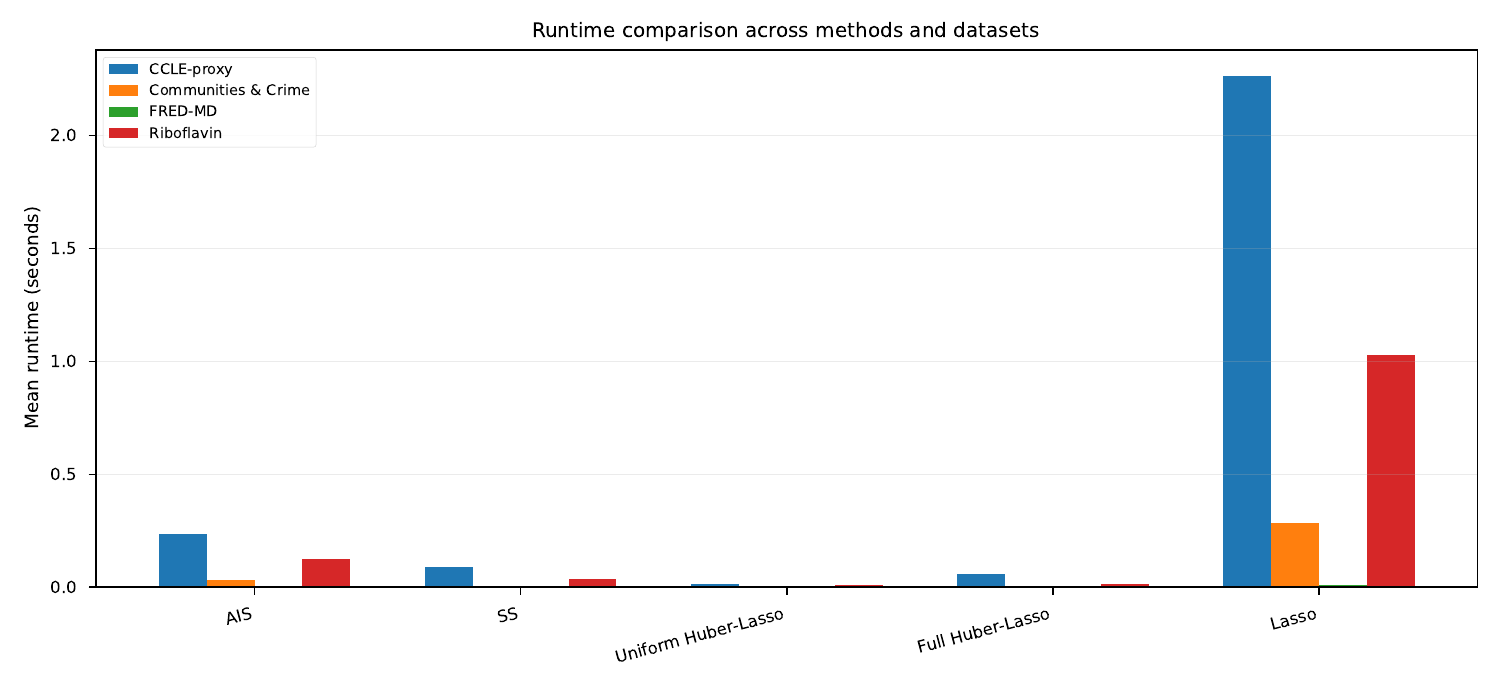}
  \caption{\emph{Left:} Log-log test MSE vs.\ $m$ on four real datasets.
    \emph{Right:} Wall-clock runtime.
    AIS is $10$--$100\times$ slower than Uniform HL per call;
    SS is the fastest method.}
  \label{fig:real}
\end{figure}

\subsection{Real Data (Table~\ref{tab:real}, Figure~\ref{fig:real})}

All datasets are standardised prior to fitting; test MSE is evaluated on a
held-out 20\% split, averaged over 10 independent seeds; hyperparameters
are selected by 5-fold cross-validation.

\paragraph{Riboflavin ($n{=}71$, $p{=}4{,}088$)(~\cite{buhlmann2014}).}
This dataset represents an extreme $p\gg n$ regime. AIS achieves a convergence
slope of $-0.793$ in log-log test MSE versus $m$, and attains 29.5\% lower MSE
than Uniform Huber-Lasso at $m{=}22$. SS collapses to a near-zero slope
($-0.063$): with only $n{=}71$ total observations, each stratum contains at most
$n_k\le5$ observations, which violates the proportional allocation requirement
$n_k\asymp n/K$ stated in Proposition~\ref{prop:ss_mom}, and causes the
geometric-median aggregation to degenerate.

\paragraph{Communities \& Crime ($n{=}319$, $p{=}122$).}
AIS achieves a convergence slope of $-0.529\approx-0.5$, closely tracking the
theoretical rate. Differences between AIS and Uniform Huber-Lasso at $m\ge76$
fall within one standard deviation.

\paragraph{CCLE-proxy ($n{=}500$, $p{=}5{,}000$, 8\% contamination).}
All methods exhibit shallow convergence slopes ($-0.044$ to $-0.021$), reflecting
the dominance of the irreducible $O(\varepsilon)$ contamination bias
from~\eqref{eq:cont} at these subsample sizes. AIS achieves the lowest test MSE
at every value of $m$.

\paragraph{FRED-MD ($n{=}399$, $p{=}125$).}
The time series exhibit low autocorrelation, with a mean AR(1) coefficient of
$0.005$. All methods reach the oracle MSE level at $m{=}63$, and the
$\alpha$-mixing correction from Theorem~\ref{thm:mixing} is negligible in
practice ($M\in[1,5]$, $B{=}12$).

\begin{table}[t]
\centering
\caption{Real-data test MSE (mean$\pm$std, 10 seeds). Bold: best subsampled.}
\label{tab:real}
\setlength\tabcolsep{2pt}
\scriptsize
\begin{tabular}{llcccc}
\toprule
Dataset & $m$ & AIS & SS & Unif HL & Full HL\\
\midrule
\multirow{4}{*}{\makecell[l]{Riboflavin\\$n{=}71,p{=}4088$}}
 & 5  & .604{\tiny$\pm$.085} & .616{\tiny$\pm$.055} & .596{\tiny$\pm$.102} & \multirow{4}{*}{.137}\\
 & 11 & \textbf{.375}{\tiny$\pm$.159} & .597{\tiny$\pm$.042} & .470{\tiny$\pm$.127} &\\
 & 16 & \textbf{.214}{\tiny$\pm$.115} & .589{\tiny$\pm$.031} & .403{\tiny$\pm$.104} &\\
 & 22 & \textbf{.201}{\tiny$\pm$.086} & .555{\tiny$\pm$.042} & .285{\tiny$\pm$.092} &\\
\midrule
\multirow{4}{*}{\makecell[l]{Comm.\&Crime\\$n{=}319,p{=}122$}}
 & 25  & .056{\tiny$\pm$.019} & .063{\tiny$\pm$.008} & \textbf{.046}{\tiny$\pm$.012} & \multirow{4}{*}{.023}\\
 & 51  & \textbf{.037}{\tiny$\pm$.007} & .049{\tiny$\pm$.007} & .035{\tiny$\pm$.013} &\\
 & 76  & .029{\tiny$\pm$.007} & .048{\tiny$\pm$.009} & .029{\tiny$\pm$.005} &\\
 & 102 & \textbf{.028}{\tiny$\pm$.004} & .038{\tiny$\pm$.004} & .027{\tiny$\pm$.003} &\\
\midrule
\multirow{4}{*}{\makecell[l]{CCLE-proxy\\$n{=}500,p{=}5000$}}
 & 40  & \textbf{54.0}{\tiny$\pm$.5} & 54.1{\tiny$\pm$.8} & 56.1{\tiny$\pm$1.8} & \multirow{4}{*}{40.4}\\
 & 80  & \textbf{53.5}{\tiny$\pm$1.1} & 53.7{\tiny$\pm$.6} & 56.4{\tiny$\pm$2.6} &\\
 & 120 & \textbf{51.2}{\tiny$\pm$.9} & 54.3{\tiny$\pm$.9} & 55.6{\tiny$\pm$4.3} &\\
 & 160 & \textbf{51.1}{\tiny$\pm$2.6} & 53.4{\tiny$\pm$.9} & 54.3{\tiny$\pm$4.8} &\\
\midrule
\multirow{4}{*}{\makecell[l]{FRED-MD\\$n{=}399,p{=}125$}}
 & 31  & \textbf{7.8e-5}{\tiny$\pm$3.0e-5} & 1.0e-4 & 8.6e-5{\tiny$\pm$2.6e-5} & \multirow{4}{*}{1.0e-4}\\
 & 63  & 1.0e-4 & 1.0e-4 & 1.0e-4 &\\
 & 95  & 1.0e-4 & 1.0e-4 & 1.0e-4 &\\
 & 127 & 1.0e-4 & 1.0e-4 & 1.0e-4 &\\
\bottomrule
\end{tabular}
\end{table}

\section{Conclusion}

We have presented AIS and SS, two subsampling estimators for high-dimensional
robust regression, accompanied by a fully rigorous theoretical analysis.
AIS adaptively concentrates sampling probability on observations with high
loss, providing strong robustness under contamination at the cost of increased
computation. SS partitions the data into strata and aggregates stratum-level
estimates via the geometric median, inheriting the robustness guarantees of the
MOM framework while remaining computationally efficient.

On the theoretical side, we have established finite-sample error bounds achieving
the minimax-optimal rate $O(\sqrt{s\log p/m})$ under sub-Gaussian design and
finite-variance noise, an explicit $O(\varepsilon)$ contamination bias bound, a
corrected $\alpha$-mixing extension using the calendar-time block protocol,
and a fully specified de-biased asymptotic normality result enabling valid
coordinate-wise confidence intervals.

\paragraph{Future directions.}
Several open problems remain.
First, a martingale stability analysis of all AIS iterates would provide
convergence guarantees for intermediate rounds of the algorithm, not just at
termination.
Second, an information-theoretic lower bound that separates AIS from uniform
subsampling under contamination would clarify the fundamental limits of adaptive
sampling in this setting.
Third, extensions to generalized linear models and nonparametric regression would
broaden the applicability of the framework.
Fourth, improved aggregation strategies for SS in the small-strata regime would
address the empirical failure mode observed on the Riboflavin dataset.
Fifth, federated learning settings, where data are distributed across multiple
nodes and communication is costly, provide a natural application domain for
stratified subsampling.

\bibliographystyle{unsrtnat}
\bibliography{references}

\end{document}